# A Reconstruction, Assessment, Error Analysis and Simulation of a Method for Measuring π That Could Have Been Used 3000 Years Ago.

By: David Neustadter Ph.D.

## Abstract

There is little known about the methods used by the ancient Babylonians and Egyptians to arrive at their recorded estimates of the value of Pi.  A surprisingly accurate estimate of Pi was recently revealed coded within a verse in the book of 1 Kings, the value of which suggests how it might have been measured.  The coded value is 111/106, which is a continued fraction representation of Pi/3.  This suggests that the value may have been measured using an iterative measurement of remainders when comparing the two lengths C (circumference of the circle) and 6R (6 times the radius).  This article describes a method that could have been used 3000 years ago to make such a measurement, the expected measurement errors, and a computer simulation that assesses the chances of such a method succeeding in obtaining the coded value of 111/106.  The result indicate that with the technology available at the time the proposed measurement method would have been possible and would have had about a 75% chance of producing the result 111/106 after a few hundred measurements.

## Introduction

A relatively recent addition to the discussion of the ancient recorded values of Pi is a surprisingly accurate estimate of Pi coded within the very Biblical verse (1 Kings, chapter 7, verse 23) that is often referenced as providing a particularly poor estimate.  The hidden value was first published by Matityahu Munk in 1962 (Munk, Three Geometric Problems in the Bible and the Talmud, 1962) and (Munk, The Halachic Way in Solving Special Geometric Problems, 1968).  The value is coded using the "gematria" coding in which numerical values are assigned to each letter in the Hebrew alphabet[1].  The coded value, the ratio of the written spelling to the pronounced spelling of the word "kav", is 111/106, which is a highly accurate approximation of Pi/3 (111/106 = 1.047169… ; Pi/3 = 1.047197…).  The value of Pi which corresponds to this coded value is Pi = 3.141509, which is significantly more accurate than any estimate of Pi known to have existed until hundreds of years after the book of 1 Kings was written[2].

---

[1] The "gematria" alphabetical numeral system is known to have been used for numbering in Hebrew since about 125 BC (Chrisomalis, 2010, p. 52), and its use in the book of 1 Kings, which is traditionally believed to have been written by Jeremiah the Prophet about 500 years earlier, is discussed later in this article.

[2] The book of 1 Kings is traditionally believed to have been written by Jeremiah the Prophet who prophesied from 627 BC to at least 586 BC when the first Jewish Temple was destroyed.

Belaga (Belaga, 1991) pointed out that this ratio is a continued fraction representation of Pi/3. This fact is what led to the idea that this value may have been arrived at using the measurement method proposed here which produces one additional level of the continued fraction with each measurement iteration. This article assesses the feasibility of the theory that the engineers who built the First Temple in Jerusalem, or some other similarly well-funded engineers of that time period, could have measured the value of 111/106 as an estimate of Pi/3.

Of course, the people making the measurement would not have needed to know about continued fractions in order to arrive at this measurement method. It is the natural method that one would use when trying to find the relationship between two lengths. In fact, the described method is an application of the Euclidean Algorithm, described by Euclid in *Elements* written around 300 BC.

To address the feasibility of this measurement method having been used 3000 years ago and having led to the ratio of 111/106 as an estimate of pi/3, I present here the reconstruction of a measurement method as it may have been performed 3000 years ago, an analysis of the different types of measurement error, and a simulation of the expected measurement results.

## Reconstruction of the measurement method

The proposed measurement method is composed of the following steps:

1) Use a stick with a sharp nail at each end as a "radius" to make a perfectly circular groove in a large flat surface of hard clay. They may have found, as I did, that using a hard metal insert with a small indentation at the center of the circle makes it easier to keep the center of the circle precisely in place as the groove is carved.
2) Use a wire stretched taught between two pegs as a guide to cut a straight groove in a large flat surface of hard clay with a length longer than 6 times the length of the "radius" stick.
3) Use the "radius" stick with nails at each end to mark off the length of 6 radii (6R) (3 times the diameter, and also the perimeter of the inscribed hexagon) along the straight groove.
4) Lay a fine wire of soft metal into the circular groove and cut the wire with a sharp knife so that it fits perfectly around the circumference of the circle. This is easily done by widening a small section of the groove and allowing the two ends of the wire to overlap, lying next to each other in the widened section, and then cutting them both simultaneously with a sharp knife. This wire now has a length C.
5) Remove the soft metal wire with length C from the circular groove and lay it into the straight groove. Cut the wire at the 6R mark to produce two pieces of wire with lengths 6R and (C-6R).
6) Using the piece of wire with length (C-6R) as a guide, cut multiple pieces of length (C-6R) and place them adjacent to each other in the groove to measure how many times (C-6R) fits into 6R and with what remainder. Since (C-6R) fits into 6R 21 times, we will call the remainder (6R – 21(C-6R)).

7) Place a piece of wire in the groove next to the 21 pieces of length (C-6R) and cut it at the 6R mark. This piece will be of length (6R – 21(C-6R)). Remove the 21 pieces of length (C-6R) from the groove, and using a piece of wire of length (C-6R) make a mark along the groove at the length (C-6R). Using the piece of wire with length (6R – 21(C-6R)) as a guide, cut multiple pieces of length (6R – 21(C-6R)) and place them adjacent to each other in the groove and compare them to the marked length of (C-6R) to measure how many times (6R – 21(C-6R)) fits into (C-6R). At this point we will ignore any remainder because, as we will see, the answer at this point will vary from measurement to measurement, so the measured remainder is obviously just random noise. Because we are ignoring the remainder we will round the measured value to the closest whole number; if less than a half of a piece is needed to reach (C-6R) we will round down and if more than a half of a piece is needed to reach (C-6R) we will round up.
8) Record the result of step 7 and repeat steps 4 through 7 many times to see what result comes out most often. The distribution of the results and their assessment are discussed below.

## The tools of measurement

The tools required for this measurement method include soft metal wire, sharp nails, and a sharp cutting blade. Evidence of the availability of both soft metal wire and sharpened hard metal tools during the early First Temple period (around 1000 BC) can be found in both the Bible itself and in archeological sources.

The Bible, in Exodus 39:3 describes how fine gold wire was made and woven into the priestly garments at the time of the exodus from Egypt (The Holy Scriptures, 1988, p. 111), about 500 years before the building of the First Temple. Archeological evidence suggests that fine soft metal wire was used for jewelry and decoration as far back as 2000 BC (Newbury & Notis, 2004). There is some debate among archeologists about how this fine wire was made, and in particular, at what point the process of wiredrawing was developed. Newbury and Notis (Newbury & Notis, 2004) describe the methods of making wire that were used before wiredrawing, and suggest that wiredrawing using soft metal dies may have been used as far back as 2000 BC to produce wires as thin as 0.5mm in diameter. However it was made, the evidence indicates that it is fair to assume that around the time of the construction of the First Temple they would have been able to produce 0.5mm wire out of copper, silver or gold.

In my experimentation with the measurement method, I found that 0.5mm copper wire was convenient for the measurement in terms of being flexible enough to easily place into both circular and straight grooves, yet rigid enough that it held its shape without needing to be pulled. 0.5mm copper wire is also easily cut with a sharp knife. Based on the properties of 0.5mm copper wire and the archeological evidence that they would have been able to produce such wire 3000 years ago, I made the error assessments and measurement simulations described below assuming that 0.5mm copper wire was used for the measurements.

Bronze tools, including sharpened tools such as needles, arrowheads and knife blades, have been found in Egypt dating back to at least the Eighteenth Dynasty, around 1200 BC (Petrie,

2013) (The Petrie Museum of Egyptian Archeology, n.d.). Iron tools, although less common, also have been found in Egypt dating from the period of the Eighteenth Dynasty (The Petrie Museum of Egyptian Archeology, n.d.). According to Scheel, iron was imported to Egypt during this period, and was produced locally only much later (Scheel, 1989). Biblical stories predating the First Temple, such as those describing the swords used by Ehud and Yoav[3], also support the existence of sharp blades among the Hebrews of the time. It therefore seems fair to assume that engineers at the time of the building of the First Temple would have been able to obtain or produce sharpened bronze or iron nails and knives for the purpose of making these measurements.

## Assessment of measurement errors

In order to assess the different measurement errors that would be expected from the proposed measurement method, I reproduced all of the steps of the measurement method and measured the error that resulted from each step.

The tools I used to assess the measurement errors included a microscope, a calibrated 1 meter long caliper with 0.02 mm precision, and a Micro-Vu microscopic measurement device (Micro-Vu, Windsor, California) with better than 0.01 mm precision.

Note: The measurement error needed to cause them to arrive at the wrong answer at the accuracy measured (to measure 4 in the second iteration instead of 5) would be an error of approximately 0.003% in the measurement of the circumference. As discussed below, the smallest radius circle they would have been likely to have used for these measurements would have been about 200mm. For this size circle, the measurement error necessary to cause them to measure the wrong answer at the accuracy measured would be about 38 microns. Any errors that have a magnitude significantly smaller than this can therefore be ignored, as they would not have been expected to affect the measurement accuracy.

The following potential sources of error were assessed:

1) The change in the shape of the cross-section of the wire when it is bent into a circle (if the cross-section changes shape, the wire will not sit centered in the groove, which would change the measured circumference of the circle).
    a. The distortion of the cross-section of a wire during bending depends on the properties of the material; in particular, how much the material stretches and/or compresses during bending.
    b. To assess this source of error, I bent a 3mm diameter copper wire into a circle of ~20mm diameter and cut through the wire to assess the cross-section shape. There was no discernable difference in cross-section shape between the straight wire and the bent wire. This potential source of error was therefore discounted.

---

[3] See Judges chapter 3 and II Samuel Chapter 20. Additionally, see I Samuel chapter 13 which tells of a period during the reign of King Saul when the Philistines had a monopoly on hard metalworking and would not allow the Hebrews to acquire swords.

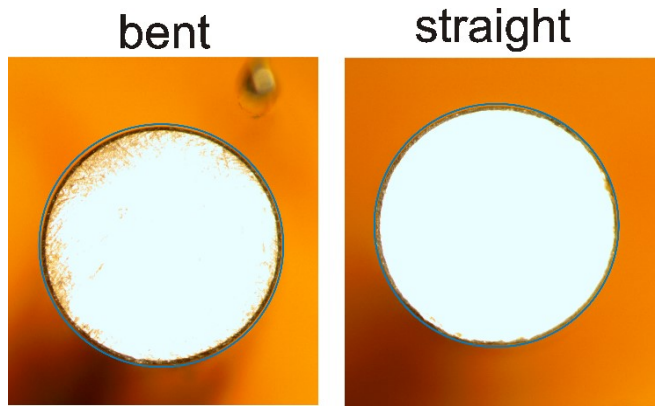

Cross sections of bent and straight wire photographed through a microscope. The overlaid blue lines are perfect circles. Note that there is no discernable distortion of the circular cross section.

2) The change in length of a piece of wire when it is bent into a circle versus when it is straight.
    a. The change in length of a wire during bending depends on the properties of the material; in particular, how much the material stretches and/or compresses during bending. A wire made of material that can stretch but can't compress will lengthen during bending and a wire made of material that can compress but can't stretch will shorten during bending.
    b. To assess this source of error, I made two marks on a 3mm diameter copper wire and then bent it into a circle of ~20mm diameter and compared the distances between the marks along the inner surface and outer surface of the circle with the distance between the marks on the straight wire before it was bent.
    c. The shortening of the inner edge of the bent wire was slightly more than the lengthening of the outer edge of the bent wire, leaving the length of the midline of the wire (which is what lies in the groove and determines the measured circumference of the circle) slightly shorter when bent than when straight.
    d. This shortening should be proportional to the diameter of the wire. The measured shortening of the midline was 171 microns for a 3mm diameter wire bent into a full circle. In the simulation, the value used for this error is therefore 57 microns per millimeter of wire diameter.

3) The change in length of a wire when it is cut using a sharp knife.
    a. In order to assess the change in shape and length of the end of a copper wire due to cutting with a sharp knife, I cut numerous wires and measured their lengths before and after cutting.

b. In addition, I performed a slow motion cut under a microscope to analyze the changes in shape and length during the cut in order to understand how the material compresses and stretches as it is cut by a knife blade.
c. The shape and length of the cut wire depends on the blade angle of the knife. Since a 20 degree blade angle is currently considered the standard optimal balance between sharpness and sturdiness for normal cutting, I assume a 20 degree blade angle for the analysis and simulation.
d. Based on the measurements described above, the simulation assumes the use of a 20 degree blade angle on a 0.5mm diameter wire which produces a wire end with a 20 degree bevel, the short side of which is 0.085mm shorter than the initial cut position, and the long side of which protrudes past the initial cut position by 0.095mm.

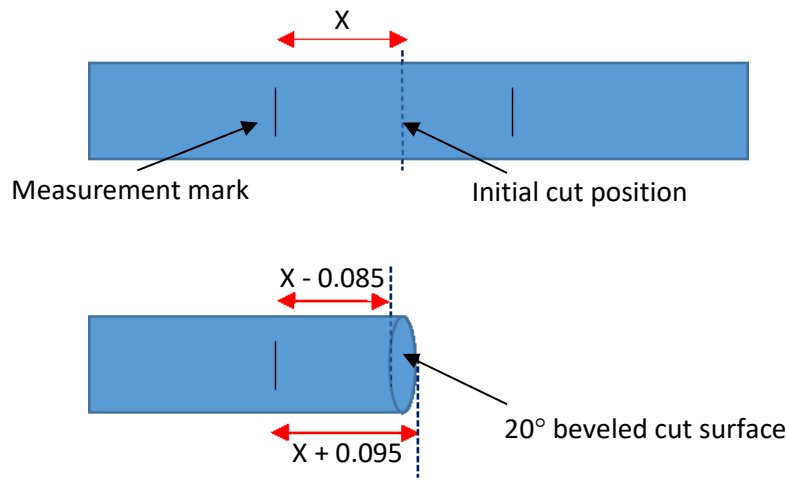

4) The accuracy with which multiple cut pieces of wire can be laid adjacent to each other in a groove.
   a. As noted above, the ends of wire cut with a sharp knife are beveled. For this reason, it is expected that there will be a uniformly distributed random error when placing such wires adjacent to one another since the amount of overlap of the beveled ends will depend on the relative orientation of the angled tip of both pieces of wire.
   b. The length errors were measured when cutting a wire and then placing the two cut pieces adjacent to each other in a groove. The mean and distribution of the errors were found to be consistent with this error being a uniformly distributed error resulting from the relative orientation of the two beveled ends of the wires. The results were consistent with the wire cutting measurements above which suggest that the beveled edges would result in a uniformly distributed elongation due to cutting and juxtaposition ranging from 0.01mm to 0.19mm.
   c. Since the simulation already includes the elongation due to cutting (0.095mm on each piece for a total of 0.19mm), the additional error due to juxtaposition is a shrinkage which is a random uniformly distributed value between 0mm and 0.18mm for every juxtaposition of two cut pieces of wire.

5) The accuracy with which a stick with 2 nails can be used to create a perfectly circular groove whose circumference is 2*pi times the distance between the points of the nails.
    a. When many wires are laid into the same groove and then measured, any inaccuracy of the groove itself should appear as a systematic error in the measurement of the circumference that repeats itself in repeated measurements.
    b. After taking into account the lengthening due to cutting discussed above, the remaining systematic error between the measured circumference and the circumference calculated from the measured distance between the points of the nails was smaller than the expected error of the mean due to the random noise between measurements. This error was therefore assumed to be a result of the random measurement noise (which is discussed below) and was not included in the simulation as a distinct error.

6) The accuracy with which a soft metal wire can be laid into a circular groove to precisely measure its circumference (the wire can be held in place within the groove using balls of clay placed periodically around the circle).
    a. When many wires are laid into the same groove and then measured, any inaccuracy resulting from the way the wire lies in the groove should appear as random error, resulting in a different measurement error for each measurement.
    b. The random error in the measured circumference, measured as the standard deviation among many measurements, was found to change with the radius of the circle. By fitting the random error relative to the radius to a straight line, the relationship between the random error and radius was estimated to be standard deviation = 0.05 + 8.68e-4 * radius. In the simulation, this was used as the formula for the standard deviation of a normal distribution of error values from which a random error was produced for each simulated circumference measurement.

7) The accuracy with which a piece of wire can be cut to match the length of another piece of wire being used as a guide.
    a. Repeated measurements were made of the accuracy with which one wire can be cut using a sharp knife to match the length of another wire.
    b. The standard deviation of the error in cutting a matching length was found to be 0.09mm. This value was used as the standard deviation of a normal distribution of error values from which a random error was produced for each simulated cut of wire to match the length of another piece of wire.

8) When a stick with 2 nails is used to make a circular groove and to mark off 6 times its length along a straight groove; with what level of accuracy is the straight length equal to 6 times the radius of the circular groove.
    a. A stick with 2 nails was used to make a circular groove and to mark off a length of 6 times the length of the stick along a straight groove.
    b. The marked length differed from 6 times the measured radius of the circular groove by less than 5 microns, an error which would have no impact on the accuracy of the proposed measurement. This potential source of error was therefore ignored in the simulation.

# The Simulation
## Measurement parameters and recording of results

In addition to the wire diameter and the different types of errors and noise in the different stages of the measurement, there are two additional parameters that must be assumed for the simulation: the radius of the circle used for the measurement and the number of times to repeat the measurement before selecting the most common result as the correct answer.

### *The radius of the circle used for the measurement*

As noted above in the discussion of the different types of measurement errors and noise, different types of errors behave differently as the size of the measurement circle changes. The table below describes the behavior of the different types of errors:

| category | Examples |
| --- | --- |
| Fixed error that does not scale with radius | 1) Change in the length of a wire when bent into a circle<br>2) Change in length of wire when cut with a knife blade |
| Random error that does not scale with radius | 1) Elongation due to juxtaposition of pieces of wire with beveled ends<br>2) Random error in cutting a piece of wire to match the length of another piece of wire |
| Fixed error that scales with radius |  |
| Random error that scales with radius | 1) Random error in placement and cutting of wire in a circular groove |

As can be seen from the table, most of the error types do not scale with the radius, which implies that as the radius increases, the measurement becomes more accurate (the error becomes smaller relative to the measured values).

It is also important to note the difference between fixed errors and random errors. As they had no accurate value of Pi to which to compare their measurements, people making measurements of Pi 3000 years ago would have had no way to identify the existence of fixed errors in their measurements.  These fixed errors, if they were large enough, would simply have caused them to arrive at the wrong answer without any way of detecting the error.

Random errors, on the other hand, would have been obvious to them.  People who made measurements on a regular basis would have been well aware that upon repeating measurements they obtain varied results, and they could easily have recorded those results and found that certain results occur more often than others.  It is reasonable to assume that when making measurements with random error they would have recorded the results of repeated measurements by counting how many times each result was obtained and forming a histogram plot using stacks of rocks, scratches on sticks, or some similar method of counting the number of results in each "bin".  Using this method they would easily identify the fact that the histogram has a distribution about a peak (the most common value) and would presumably select the most common value as the correct result.

When setting out to measure Pi according to the proposed method, they would presumably have started out using a small circle to make the measurement faster and easier.  As such, they would have found that already on the first iteration of the measurement when measuring how many times (C-6R) fits into 6R (in step 6 in the method described above), they would have obtained varied results due to measurement noise.  By trying different sized circles, they would have found that the percentage of the measurement results that were not equal to 21 changed consistently with the size of the circle.  Since they wanted to proceed to the second iteration (step 7 of the method described above) in order to obtain higher accuracy, they would have wanted to use a circle large enough to provide a result of 21 on the first iteration close to 100% of the time.

The plot below shows the results of simulated measurements using the measurement error values described above.  As can be seen in the plot, the result of the first iteration of the measurement becomes more consistent as the radius increases.  The plot indicates that at a radius of about 450mm the correct result is obtained about 98% of the time and the graph levels off, increasing very slowly for additional increase in radius.  For this reason, 450mm was used as the radius of the measurement circle in the simulation (sensitivity of the simulation results to the chosen radius was analyzed and is discussed below).

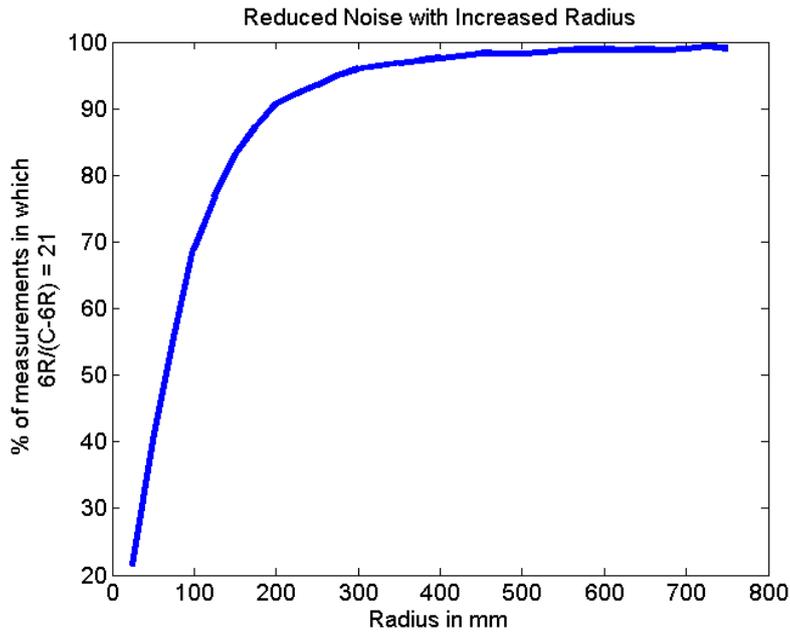

### The number of times to repeat the measurement

Since this measurement method will result in noisy results for the second iteration, it is necessary to assume, as described above, that the people making the measurement were familiar with the concept of noisy measurements. Presumably they would record the results in some sort of histogram plot and then choose the most common result. The question is; how would they have decided how many measurements to make before selecting the most common value as the result? Clearly, the more measurements they make, the more likely it will be that the most common value is indeed the correct result, but how would they have known when they had made sufficient measurements to be confident that they will get the correct result?

Presumably, they were not able to analyze the statistics of the distribution and assess the probability of getting the correct answer after a given number of measurements. However, they were able to graphically visualize the plot of the distribution of the results, and they could have been aware of the smoothness of the plot. They may, for instance, have known from experience that the plot should have a single peak and decrease monotonically and smoothly on both sides of that peak. It would be reasonable to suggest that they may have continued making measurements until they obtained a plot in which all values above some baseline noise were monotonically decreasing away from a peak value. The simulation assumes this approach and repeats the measurement procedure until a histogram meeting this criterion is obtained. (Assessment of the likelihood of obtaining the correct answer with fewer or more measurements was also assessed and is presented below.)

This plot shows a typical result of the second measurement iteration (step 7) from 300 measurements, from which the correct answer of 5 would be easily identified:

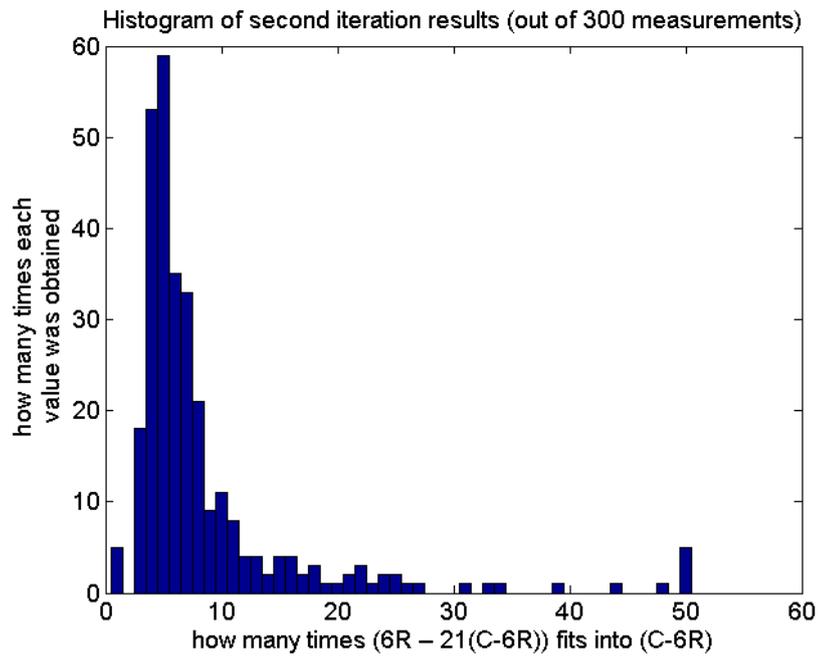

The same histogram (with bins from 1 through 16, right to left, and numbered using Hieratic Numerals) is shown as a demonstration of how such a plot may have looked 3000 years ago scratched into a slab of clay.

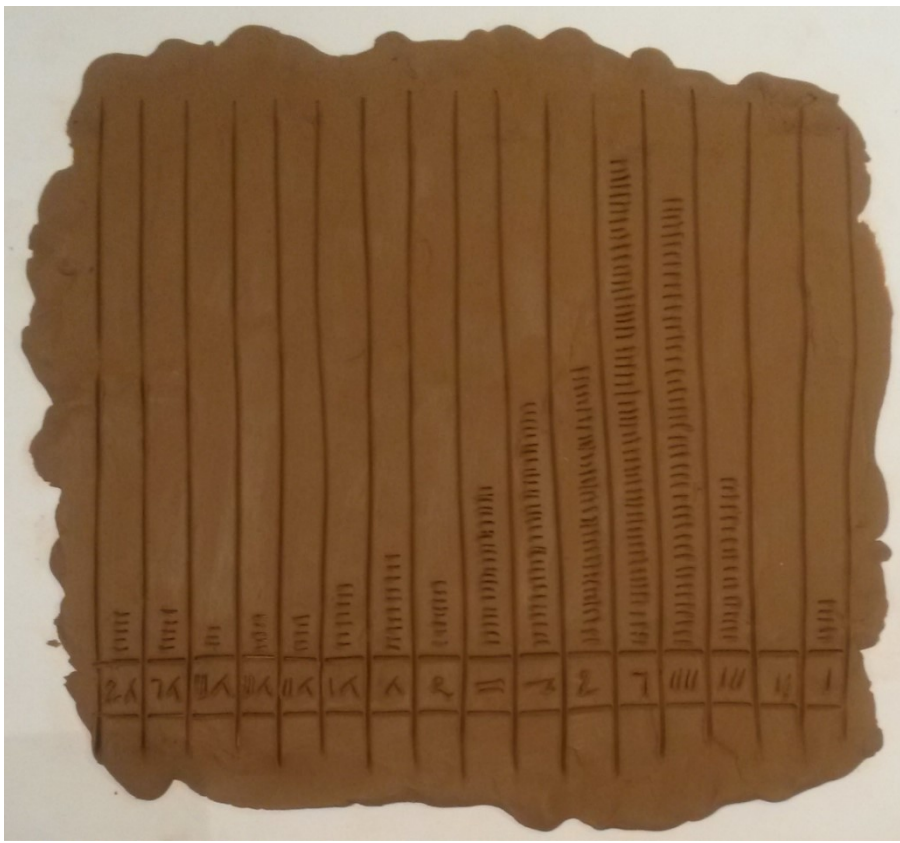

What now remains to be assessed using the simulation is; what are the chances that, given the measurement errors described above, a wire diameter of 0.5mm and a circle with a radius of 450mm, they would have obtained the result of 5 in the second measurement iteration often enough to lead them to conclude that it was the correct value. This will depend on two factors; will the random noise be small enough to get consistent enough results to be able to recognize the peak of the distribution, and will the fixed errors be small enough for the peak of the distribution of the second iteration of the measurement to be correct. This was assessed using the simulation described below.

## The simulation algorithm

The simulation is written in MATLAB (MathWorks, Natick, MA) and simulates the measurement method as it could have been performed 3000 years ago by performing the following steps:

1) The variable c_minus_six_r_ref, representing the length of the piece of wire obtained by cutting the circumference wire at the length of 6 radii, as in step 5 of the measurement method described above, is calculated by adding together the following values:
    a. 2 * pi * 450 (the circumference of the circle) minus 6 * 450 (6R)
    b. The wire bending error (0.057 * wire diameter = 0.0285)
    c. The circumference length random error (0.3538) times a random value from a normal distribution with mean=0 and stdev=1
    d. 3 times the cut elongation error (0.285) (2 for the two ends of the wire cut to be the circumference of the circle, and a third for the cut made at the length of 6R)
    e. The cut length random error (0.09) times a random value from a normal distribution with mean=0 and stdev=1
2) Step 6 of the measurement method, measuring how many times (C-6R) fits into 6R, is simulated by accumulating "pieces of wire" cut to match the original wire of length c_minus_six_r_ref until the total length is greater than 6R and then subtracting one. The accumulated length is constructed as follows:
    a. Initial length = c_minus_six_r_ref
    b. Additional "pieces" are added by adding the following values and repeating until the total length is greater than 6R:
        i. c_minus_six_r_ref
        ii. Cut length random error (0.09) times a random value from a normal distribution with mean=0 and stdev=1
        iii. Juxtaposition error which is a random value from a uniform distribution between -0.18 and 0
    c. The number of times (C-6R) fits into 6R is recorded as one less than the number of "pieces" required to pass the length of 6R
3) new_piece_ref, the remainder from step 6, (6R – 21(C-6R)), which according to the measurement method is obtained in step 7 by cutting the last (C-6R) piece at the 6 radius mark, is simulated by adding together the following values:

a. 6 * 450 minus the accumulated length of the pieces in the previous step not including the last piece that made it longer than 6R
b. Juxtaposition error which is a random value from a uniform distribution between 0 and 0.18
c. Cut length random error (0.09) times a random value from a normal distribution with mean=0 and stdev=1
d. Cut elongation error (0.095)

4) The measurement of the number of times (6R – 21(C-6R)) fits into (C-6R) is then simulated by accumulating pieces of length (6R – 21(C-6R)) until their total length is greater than (C-6R), and then subtracting one if less than half of the last piece extends past the (C-6R) mark (rounding off to the nearest integer). The accumulation is simulated as follows:
   a. Initial length = new_piece_ref
   b. Additional "pieces" are added by adding the following values until the total length is greater than (C-6R):
      i. new_piece_ref
      ii. Cut length random error (0.09) times a random value from a normal distribution with mean=0 and stdev=1
      iii. Juxtaposition error which is a random value from a uniform distribution between -0.18 and 0
   c. If less than half of the last piece extends past (C-6R) then the number of times (6R – 21(C-6R)) fits into (C-6R) is recorded as the number of "pieces" required to pass the length of (C-6R) minus one, otherwise the number of times (6R – 21(C-6R)) fits into (C-6R) is recorded as the number of "pieces" required to pass the length of (C-6R).

5) The simulated measurement is repeated and the results recorded until the distribution of the results meets the predetermined stopping criteria, at which point the peak value is selected as the correct value. The stopping criteria were defined, as discussed above, to indicate a single identifiable peak and a distribution that is monotonically decreasing on both sides of the peak:
   a. The peak must be at least 5 counts and at least 5% larger than the values of the adjacent bins both below and above the peak
   b. There must be at least 5 consecutive bins, including the peak, with values above 20% of the value of the peak
   c. Within the group of consecutive bins whose values are above 20% of the value of the peak, the values must be monotonically decreasing both below and above the peak

6) The simulation can be run in a loop many times in order to assess the chances of obtaining the correct answer using this measurement method, measurement parameters and noise values.

## The Results

The results of the simulation are that the distribution of the results meets the stopping criteria after an average of about 325 measurements, and the distribution indicates the value of 111/106 (the correct value of 5 is obtained in the second iteration of the measurement) approximately 75% of the time.

These results support the plausibility that this method could have been used and would have had a 75% chance of resulting in the value of 111/106.

Even if we assume that they would not have had the patience to make 325 measurements, and would have used more relaxed stopping criteria, the chances of getting the correct results for lower numbers of measurements are still reasonable; >40% for 50 measurements and >50% for 100 measurements.

# Interesting Features of the Measurement Method and Results
## Effects of Fixed Errors and Random Errors

An interesting feature of this measurement method relates to the fact that it is influenced differently by different types of measurement errors. Of particular interest is the effect of random errors. Ordinarily, random errors are expected to average out with large numbers of iterations and not to have any effect on the mean of the distribution of the results. However, in this measurement method that is not the case.

The measurement that we are making is the number of times one length piece of wire fits into another length piece of wire. Since both pieces of wire include a random error component in their length, we are measuring the ratio of two random variables. The distribution of the ratio of two normally distributed random variables is referred to as the Reciprocal Normal Distribution. One of the properties of this distribution is that the location of both the mean and the peak of the distribution depend on the variance of the random variable in the denominator. As the variance of the distribution in the denominator increases, the mean of the Reciprocal Distribution increases and the location of its peak decreases.

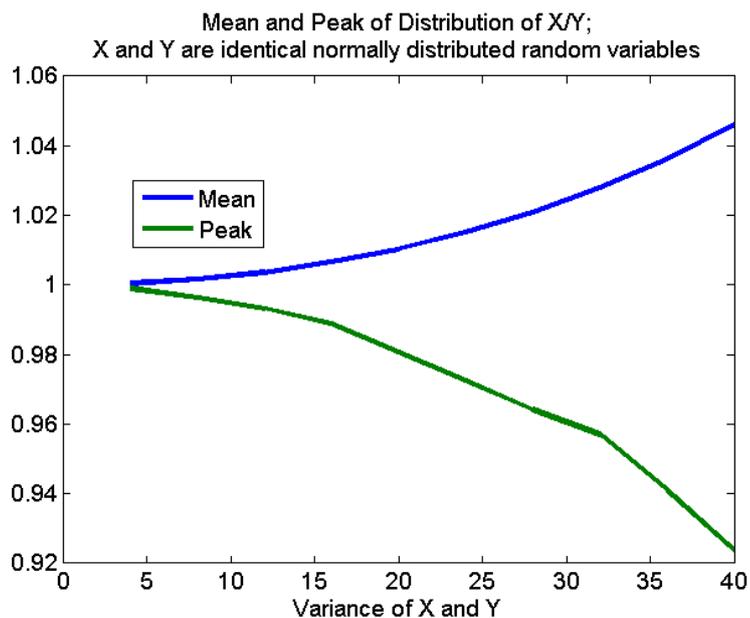

(Note: For distributions with large variance, the reciprocal normal distribution is extremely noisy. In order to obtain this smooth curve of the peak location for different variances, 10,000,000,000 random values were produced.)

For this reason, the random error in the lengths of the pieces of wire reduces the value of the peak of the distribution, inducing an error in the identification of the peak of the distribution of the results in the direction of selecting a lower peak (4 rather than the correct answer of 5).

The most significant fixed error in the measurement is the lengthening of every cut piece due to the angled wire ends resulting from being cut with a blade. As a result, each piece of length (C-6R) will be a bit longer than it should be. Because the length of (6R – 21(C-6R)) is obtained by taking the remaining part of 6R after 21 pieces of (C-6R) have been removed from it, the length of (6R – 21(C-6R)) will be smaller than it should be. The second iteration of the measurement (measurement method step 7) will therefore be assessing how many pieces of (6R – 21(C-6R)), each of which is shorter than it should be, fit into (C-6R), which is longer than it should be. This fixed error will therefore induce an error in the direction of selecting a higher peak (7 rather than the correct answer of 5).

By running the simulation without the fixed errors and again without the random errors, it was shown that indeed with only fixed errors and no random errors the result of step 7 is 7 instead of 5 100% of the time, and with only random errors and no fixed errors the result is 4 instead of 5 65% of the time. Interestingly, the values of the fixed and random errors as they were measured in this analysis cause the effects of these two errors to cancel each other out sufficiently that the correct result of 5 is obtained most of the time as described above.

## Sensitivity to Choice of Radius and Number of Measurements

The choice of radius was discussed above, and was assumed to be made based on the results of the first iteration of the measurement (measurement step 6). However, it is still interesting to investigate the sensitivity of the measurement results to the selection of radius.

Running the simulation with different radii and analyzing how the radius affects the chances of obtaining the correct result in the second iteration (measurement step 7) produced surprising results: Increasing the radius does increase the chances of obtaining the correct result, but not nearly as much as increasing the number of measurements. Investigation of this result indicated that increasing the radius has only a small effect on the results because the effects of the fixed and random errors cancel each other out, as described above, and most of the components of both of these errors scale linearly with the radius. (Note: the effect of the radius on the second iteration of the measurement was assessed assuming that the measurement only continued to the second iteration if the known correct answer of 21 was obtained in the first iteration)

What was found to have a significant effect on the chances of obtaining the correct result was the number of measurements, as discussed above.

Shown below is a plot indicating the chance of obtaining the correct result (5 in the second measurement iteration, measurement step 7) versus the number of measurements over a range of radii.

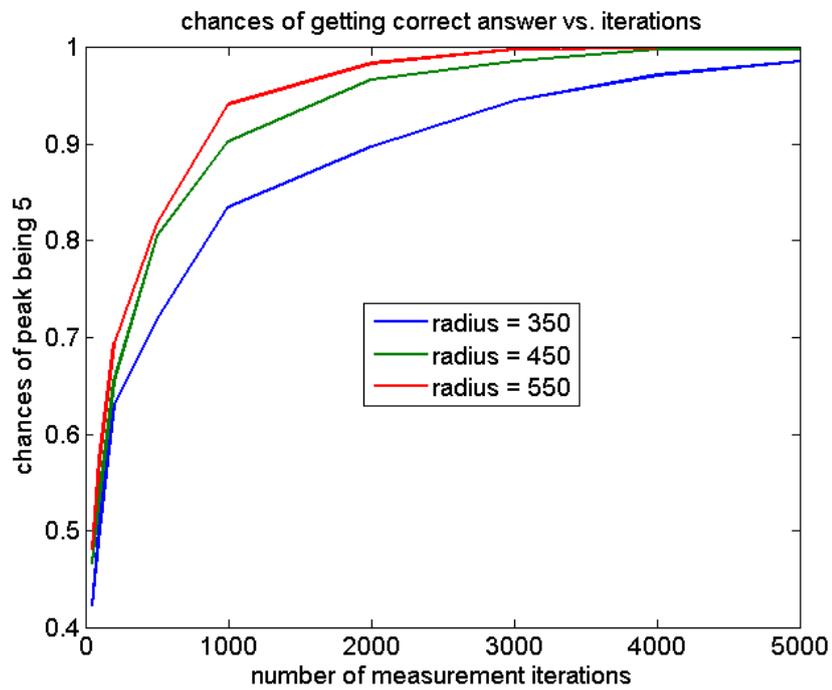